\documentclass[twoside,openany]{book}
\usepackage{amsmath,amssymb}
\newtheorem{theorem}{Theorem}[section]
\newtheorem{definition}[theorem]{Definition}

\def \M {\mathcal{M}}

\def\al{\aligned}
\def\eal{\endaligned}
\def\be{\begin{equation}}
\def\ee{\end{equation}}
\def\befr{\begin{frame}}
\def\enfr{\end{frame}}
\def\lab{\label}
\def\a{\alpha}

\def\n{\nabla}

\def\R{{\bf R}}
\def\M{{\bf M}}

\def\al{\aligned}

\def\p{\partial}
\def\d{\nabla}

\newpage
\vspace*{11mm}
\begin{document}

\begin{center}{\LARGE\bf Log gradient estimates for heat type equations on manifolds
}\end{center}

\markboth{\small {\sc Qi S. Zhang
}}{\small log gradient estimates}

\begin{center}
{\sc Qi S. Zhang}
\end{center}

\noindent {\small Department of Mathematics, University of California, Riverside, CA 92521, USA\\
qizhang@math.ucr.edu}

%{\small qizhang@math.ucr.edu}
\date{2023/10}

\begin{quotation}
\noindent {\bf Abstract}\quad   {\small  In this short survey paper, we first recall the log gradient estimates for the heat equation on manifolds by Li-Yau, R. Hamilton and later by Perelman in conjunction with the Ricci flow. Then we will discuss some of their applications and  extensions focusing on sharp constants and improved curvature conditions.
}

\end{quotation}

\newpage

\section{Introduction}

Gradient estimates for solutions of differential equations can be traced back at least to the study of analytic functions whose gradients in a smooth,
bounded domain can be determined by their boundary values via the Cauchy integral formula. In the real variable case, Bernstein \cite{Ber} discovered gradient bounds in a similar spirit for solutions of several elliptic equations in the interior of a domain. His method is to apply the maximum principle on a carefully chosen auxiliary function involving the solution and the modulus of the gradient, which satisfies a  differential inequality whose nonlinear term has a good sign. This general strategy was very fruitful and is still being used nowadays for wide range of equations and problems, some of which will be touched below. When the background space is a Riemannian manifold, the study of harmonic forms also involve gradient estimates in the $L^2$ and other integral sense. Hodge theory is one of the ramifications. Comparing with the Euclidean case, various curvature terms come into the play through commutation formulas such as Bochner's. This not only makes the situation more complicated but, more importantly, also reveals the connection with geometric and topological properties of the manifolds. Further applications of gradient estimates and the Bernstein technique in conjunction with commutation formulas are found in complex geometry such as Kodaira vanishing theorem and fully nonlinear elliptic equations which are beyond the scope of this paper.

\section{Li-Yau and Hamilton estimates in fixed metric case}

When the solution of an equation is positive. One can apply Bernstein's method to study the log gradient bound of the solution. This was carried out by S. T. Yau \cite{Ya} and Cheng and Yau \cite{CY} for positive harmonic functions on Riemannian manifolds. The log gradient bound is stronger than just the gradient bound. For example log gradient bound for positive harmonic functions implies the Harnack inequality which roughly states that in a compact sub-domain bounded away from the boundary, the maximum value is bounded by the minimum value times a constant which is independent of the solution. See also \cite{Mo} by Modica for a log gradient bound for nonlinear elliptic equations.

Log gradient bounds have been extended to positive solutions of the heat equation, which is the main topic of the present paper.
Let $(\M^n, g_{ij})$ be a complete Riemannian manifold. In \cite{LY:1}, P. Li an S.T. Yau discovered the following celebrated Li-Yau bound.

\begin{theorem}
Let $\M$ be a complete Riemannian $n$-manifold. Suppose  $Ric\geq-K$, where $K\geq0$ and $Ric$ is the Ricci curvature of $\M$. Let $u$ be a  positive solutions of the heat equation in $\M \times (0, \infty)$:
\be
\lab{heatequation}
\frac{\partial u}{\partial t}=\Delta u.
\ee  Then
\begin{equation}\label{Li-Yau>1}
\frac{|\nabla u|^2}{u^2} - \alpha\frac{u_t}{u} \leq \frac{n\alpha^2K}{2(\alpha-1)}+\frac{n\alpha^2}{2t},\quad \forall \alpha>1.
\end{equation}
In the special case where $Ric\geq 0$, then
\begin{equation}\label{Li-Yau.op}
\frac{|\nabla u|^2}{u^2} - \frac{u_t}{u} \leq \frac{n}{2t}.
\end{equation}
\end{theorem}

In the same paper, many applications of \eqref{Li-Yau>1} and \eqref{Li-Yau.op} have also been demonstrated by the authors, some of which are listed below.

 $\bullet$ The classical parabolic Harnack inequality,  Gaussian estimates of the heat kernel. We only state here the results when the Ricci curvature is nonnegative, which are sharp. General results can be found in the original paper by Li-Yau and the book \cite{Lib}.

 \begin{theorem} Let $\M$ be a complete noncompact $n$ dimensional Riemmannian manifold of nonnegative Ricci curvature. Let $u=u(x, t)$ be a positive solution of the heat equation $\Delta u - \partial_t u=0$, $t>0$. Then
 \[
 u(x, t_1) \le u(y, t_2)  (t_2/t_1)^{n/2} \exp (d^2(x, y)/(4(t_2-t_1)), \quad x, y \in \M, t_2>t_1>0.
 \]
 \end{theorem}

 \begin{theorem} Let $\M$ be a complete noncompact Riemmannian manifold of nonnegative Ricci curvature. Let $G=G(x, t, y)$ be the (minimal) heat kernel, i.e. $\Delta_x G - \partial_t G=0$, $t>0$ and $G(x, 0, y) =\delta(x, y)$, the Dirac delta function. For any small $\epsilon>0$, there exists a positive constant $C$, depending only on $\epsilon$ and the dimension such that
 \[
 \frac{1}{C |B(x, \sqrt{t})|} \exp(- \frac{d^2(x, y)}{4(1-\epsilon) t}) \le G(x, t, y) \le \frac{C}{|B(x, \sqrt{t})|} \exp(- \frac{d^2(x, y)}{4(1+\epsilon) t})
 \]Here $d=d(x, y)$ is the geodesic distance between $x$ and $y$.

 \end{theorem}

 $\bullet$  Estimates of eigenvalues of the Laplace operator, and estimates of the Green's function of the Laplacian.

$\bullet$ Moreover
\eqref{Li-Yau.op} coupled with the small time asymptotic expansion of the heat kernel can even imply the Laplacian/volume Comparison Theorem (see e.g. \cite{Chowetc} page 394).

$\bullet$ Later P. Li \cite{Li} himself used the Harnack inequality to study the fundamental group of noncompact manifolds. He proved the following theorem.

\begin{theorem} Let $\M$ be a complete noncompact Riemmannian manifold of nonnegative Ricci curvature. Suppose the volume of geodesic balls centered at a point are maximum in the sense that $|B(p, r)|/r^n \ge const.>0$ for all $r>0$. Then the fundamental group $\pi_1(\M)$ is finite.
\end{theorem}

   This confirms a special case of Milnor's conjecture: {\it The fundamental group of complete noncompact Riemmannian manifold of nonnegative Ricci curvature is finitely generated.} See also Anderson \cite{An} for a different proof. Very recently  Elia Brue, Aaron Naber, Daniele Semola posted a paper \cite{BNS} to give a counter example. But this makes Li's result even more interesting and one can expect further generalization on the positive side of the result. Earlier positive results can be found in C. Sormani \cite{So} (linear volume growth), G. Liu \cite{Liu} (3 dimension), see also J.Y. Pan \cite{Pa} (3 dimension).

$\bullet$ See also Mabuchi \cite{Ma},  Cao-Tian-Zhu \cite{CTZ} for applications on complex geometry. These authors used Li-Yau inequality for certain (weighted) heat equation  to prove lower bound of the Green's function of certain perturbed Laplacian  which is used to study K\"ahler Einstein metrics and K\"ahler Ricci solitons. Further, some application is found by Colding and Naber in \cite{CoNa} for extending the Cheeger-Colding theory on limits of manifolds with Ricci curvature bounded from below.

We should mention that in the Euclidean space, the bound \eqref{Li-Yau.op} was obtained by Aronson-Benilan \cite{ArBe} earlier even for some porous medium equations.   Hamilton \cite{Ha93} further showed a matrix Li-Yau bound for the heat equation. This is the precursor to his matrix and trace Harnack inequality for 3 dimensional Ricci flows with positive curvature operator. Similar matrix Li-Yau bound was subsequently obtained by Cao-Ni \cite{CaNi} on K\"ahler manifolds. We will discuss these some of these results a little later.

  For the past three decades, many Li-Yau type bounds have been proved not only for
the heat equation, but more generally, for other linear and semi-linear parabolic equations
on manifolds with or without weights.  Let us mention the result by  Bakry and Ledoux \cite{BL} who
derived the Li-Yau bound for weighted manifolds by an ordinary differential inequality
involving the entropy and energy of the backward heat equation. Also P. Li-G. Tian \cite{LT} derived a Li-Yau bound for certain algebraic varieties, J.Y. Li \cite{Lij} proved mixed gradient estimates for the heat kernel.  Davies \cite{Dav}, Garofalo-Modino \cite{GM},
J.F. Li and X. Xu \cite{LX}, Qian-Zhang-Zhu \cite{QZZ}, F.Y. Wang\cite{Wan}, Jiaping Wang\cite{WanJ} and the latest Bakry-Bolley-Gentil\cite{BBG} for improved coefficients and its references.

Next we turn to another important log gradient bound  due to R. Hamilton 1993 \cite{Ha93}.

\begin{theorem}
\lab{thhament}
Let $u=u(x, t)$ be a positive solution to the heat equation $\Delta u - \partial_t u=0$ on a compact manifold $\M$ whose Ricci curvature is bounded from below by a constant $-K \le 0$. If $u \le A$, a constant, then
\[
t \frac{| \nabla u|^2}{u} \le (1+ 2 K t) u \ln \frac{A}{u}.
\]
\end{theorem}

One reason for its importance is that it is a gradient estimate related to the "density" of the entropy $-u \ln u$.
It is related to the second law of thermodynamics on the monotonicity of the entropy $-\int u \ln u \, dg$ and actually implies that the time derivative of the entropy is also monotone. So the entropy is concave in time if $Ricci \ge 0$.

Indeed, write
\be
\lab{e(t)}
E(t)= -\int_{\M} u \ln u \, dg,
\ee then using
the basic identities
\[
(\Delta - \partial_t) (u \ln u) = \frac{|\nabla u|^2}{u},
\]
\[
(\Delta  - \p_t) \frac{ |\d u|^2}{u} = \frac{2}{u} \left| u_{ij} - \frac{u_i u_j}{u}
\right|^2 + 2 R_{ij} \frac{u_i u_j}{u},
\]one deduces
\be
\lab{e'>0}
E'(t) = \int_{\M} \frac{|\nabla u|^2}{u} \, dg \ge 0,
\ee
\be
\lab{e''<0}
E''(t) = -2 \int_{\M} \frac{1}{u} \left| Hess \, u - \frac{\nabla u \otimes \nabla u}{u}\right|^2 dg - 2 \int_{\M} \frac{Ric(\nabla u, \nabla u)}{u} dg \le 0.
\ee
The proof of Theorem \ref{thhament} is to use a linear combination of the basic identities and apply the maximum principle.

Another important application of Theorem \ref{thhament} is that it allows comparison of values of the solution $u$ at different spatial points in the same time level. This is not possible from the classical Harnack inequality. See further extensions in \cite{ChH:1}.

There is a localized version of Theorem \ref{thhament}, which has a square on top of the log term. It looks worse than Hamilton's global result but it is actually necessary.

\begin{theorem}  (\cite{SZ:1})
\lab{THSZ1}
Let ${\bf M}$ be a Riemannian manifold of dimension $n \ge 2$ such that
\[
Ricci ({\bf M}) \ge -k, \qquad k \ge 0.
\]Suppose $u$ is any positive
solution to the heat equation in $Q_{R, T} \equiv B(x_0, R) \times
[t_0-T, t_0] \subset {\bf M} \times (-\infty, \infty)$.
  Suppose also $u \le A$ in  $Q_{R, T}$. Then there exists a dimensional
  constant $c$ such that
\be
\lab{1.4}
  \frac{|\nabla u(x, t) |}{u(x, t)} \le c (\frac{1}{R} +
\frac{1}{T^{1/2}}+\sqrt{k})
   \big{(} 1+ \ln \frac{A}{u(x, t)} \big{)}
\ee
in $Q_{R/2, T/2}$.
\end{theorem}

Parabolic Liouville theorem for ancient solutions of the heat equation is one of the applications of this localized gradient estimate. See \cite{LinZ}, \cite{Mos} and \cite{HY} for further development and refinement in this direction.

\section{Perelman type estimates in Ricci flow case}

The study of Li-Yau bound for heat type equations under the Ricci flow was initiated by Hamilton. In \cite{Ham1988}, he obtained a Li-Yau bound for the scalar curvature along the Ricci flow on 2-sphere. This result was later improved by Chow \cite{Chow}. In higher dimensions, both matrix and trace Li-Yau bounds for curvature tensors, also known as Li-Yau-Hamilton inequalities, were obtained by Hamilton \cite{Ham1993jdg} for the Ricci flow with bounded curvature and nonnegative curvature operator. These estimates played a crucial role in the study of singularity formations of the Ricci flow on three-manifolds and solution to the Poincar\'e conjecture. We remark that
Brendle \cite{Bre} has generalized Li-Yau-Hamilton inequalities under weaker curvature assumptions. The Li-Yau-Hamilton inequality for the K\"ahler-Ricci flow with nonnegative holomorphic bisectional curvature was obtained by H.-D. Cao \cite{Cao}.

On the other hand, the log gradient bounds in Section 2 have been extended to situations with moving metrics.
Let $g_{ij}(t)$, $t\in[0, T]$, be a family of Riemannian metrics on $\M$ which solves the Ricci flow:
\begin{equation}\label{RF1}
\frac{\partial }{\partial t}g_{ij}(t)=-2R_{ij}(t),
\end{equation}
where $R_{ij}(t)$ is the Ricci curvature tensor of $g_{ij}(t)$. One may still consider linear and semi-linear parabolic equations under the Ricci flow in the sense that in the heat operator $\frac{\partial}{\partial t}-\Delta$, we have $\Delta=\Delta_t$ which is the Laplace operator with respect to the metric $g_{ij}(t)$ at time $t$.

 The two most prominent examples are the heat equation
\be
\label{ricciheateq}
  (\Delta-\frac{\partial }{\partial t})u= 0, \, \partial_t g_{ij}= - 2 R_{ij}
  \ee and the conjugate heat equation
 \be
\label{ricciconheateq}
 H^* u \equiv (\Delta-R + \frac{\partial }{\partial t})u= 0, \quad \partial_t g_{ij}= - 2 R_{ij}.
  \ee

In 2002, in the breakthrough paper \cite{P:1}, Perelman established the equivalent of \eqref{e'>0} and \eqref{e''<0} for
\eqref{ricciconheateq}. More precisely,

\begin{definition} (F entropy \index{Perelman F entropy} and W entropy
\index{Perelman W entropy})
\lab{defFW}
Let $u$ be a positive solution to \eqref{ricciconheateq} in a compact $n$-manifold.
The F entropy is the integration of $H^*(u \ln
u)$, i.e.
\be \lab{Fshang} {\bf F} = \int_{\bf M} (\frac{|\nabla
u|^2}{u} + R u) d g(t). \ee Here $R$ is the scalar curvature.
 The W entropy is a combination of
the F entropy and the Bolzmann entropy $\bold{B} = \int u \ln u dg(t)$ together with certain scaling
factor.  Let $\tau>0$ be such that $\frac{d\tau}{dt} = -1$, define \be
\lab{Wshang} {\bf W}=\tau {\bf F} - {\bf B} - \frac{n}{2} (\ln 4 \pi
\tau) -n. \ee  i.e.
\[
{\bf W}= \int_{\bf M}  \left[ \tau (\frac{|\nabla u|^2}{u} + R u) - u \ln u
-\frac{n}{2} (\ln 4 \pi \tau) \ u  -n u \right]
d g(t).
\]
\end{definition}

\begin{theorem} (\cite{P:1})
\lab{thFMdandiao}
Perelman's $F$ and $W$ entropy are nondecreasing in time $t$. Moreover
 \[
\frac{d}{dt} {\bf F}= 2 \int_{\bf M} | Ric - Hess (\ln u)|^2 \ u \
 dg(t),
 \]
\[
\frac{d}{dt} {\bf W}=
2 \tau \int_{\bf M} \left[ Ric -Hess(\ln u) - \frac{1}{2 \tau} g \right]^2
\ u \  d g(t).
\]
\end{theorem}

The amazing thing is that $F$ and $W$ are monotone in time without any curvature condition. The monotonicity is strict unless the Ricci flow is gradient soliton.
The proof of the theorem can be based on the following identities
\[
\al
&H^*(u \ln u) = \frac{|\nabla u|^2}{u} + R u,\\
&H^*(\frac{|\nabla u|^2}{u} + R u) \\
&\quad = \frac{2}{u} \left( u_{ij}- \frac{u_i u_j}{u} \right)^2 +
 4 \nabla R \nabla u + \frac{4}{u} Ric(\nabla u, \nabla u) + 2 |Ric|^2 u  +
2 u \Delta R.
\eal
\]All curvatures terms are integrated out or enter some square terms.

In \cite{P:1}, Perelman also showed a combination of Li-Yau and Hamilton type bound for $u=u(x, \tau; y, 0)$ the fundamental solution of the conjugate heat equation \eqref{ricciconheateq} with a pole at $(y, 0)$ under the Ricci flow on a compact manifold: for $\tau=T-t$ with $T$ a fixed number,
\[
\tau \left(-2 \frac{\Delta u}{u} + \frac{|\nabla u|^2}{u^2} + R \right ) - \ln u -\frac{n}{2}
(\ln 4\pi\tau) -n \le 0.
\]

This inequality implies:
\begin{itemize}

\item Perelman's monotonicity formula for the W entropy and,

\item  his $\kappa$ non-collapsing property for volume of geodesic balls, (on his way to the Poincar\'e conjecture);

\item and after a little work,  the non-inflating property for volume of geodesic balls
(\cite{Z12},  \cite{CW13}).
\end{itemize}
These are some of the  basic properties of Ricci flows by now .

The last two statements say that if the scalar curvature is bounded in a suitable region than
\[
K_1 \le |B(x, r, t)|/r^n \le K_2,  \qquad  0 < r \le r_0,  \, 0<t \le t_0,
\]for some positive constants $K_1, K_2$ which may depend on $t_0$ and the initial metric. Here $|B(x, r, t)|$ is the volume of the geodesic ball of radius $r$ centered at $x$ with respect to the metric at time $t$. $r_0$ is a fixed scale. Note that no Ricci curvature bound is needed.

Afterwards, there have been many results on Li-Yau bounds for positive solutions of the heat or conjugate heat equations under the Ricci flow.
 \begin{itemize}

 \item For example authors of Kuang-Zhang \cite{KuZh} and X.D. Cao \cite{Cx} proved (non-sharp) Li-Yau type bound
for all positive solutions of the conjugate heat equation without any curvature condition, just like
Perelman's aforementioned result for the fundamental solution.

\item
In X.D.Cao-Hamilton \cite{CH:1} and  Bailesteanu-Cao-Pulemotov,\cite{BCP} the authors  proved various Li-Yau type bounds for positive solutions of \eqref{ricciheateq}
under either positivity condition of the curvature tensor or boundedness of the Ricci curvature.

\end{itemize}

So there is a marked difference between these results on the conjugate
heat equation and the heat equation in the curvature conditions. In view of the absence of curvature
condition for the conjugate heat equation, one would hope that the curvature conditions for the heat equation can be weakened.

 A few years ago, in Bamler-Zhang \cite{BZ}, the authors proved the following gradient estimate for bounded positive solutions $u$ of the heat equation \eqref{ricciheateq},
\begin{equation}\label{BZ}
|\Delta u|+\frac{|\nabla u|^2}{u}-aR\leq \frac{Ba}{t},
\end{equation}
where $R=R(x,t)$ is the scalar curvature of the manifold at time $t$, and $B$ is a constant and $a$ is an upper bound of $u$ on $\M\times[0,T]$. Although this result requires no curvature condition and it has some other applications such as distance distortion estimate under Ricci flows,  it is not a
Li-Yau type bound.
It turns out that a Li-Yau type bound is valid for the forward heat equation \eqref{ricciheateq}
 under the very weak assumption that the scalar curvature is bounded.

\begin{theorem}(\cite{ZZ1})\label{thmLYRF}
Let $\M$ be a compact $n$ dimensional Riemmannian manifold, and $g_{ij}(t)$, $t\in [0,T)$, a solution of the Ricci flow \eqref{RF1} on $\M$. Denote by $R$ the scalar curvature of $\M$ at $t$, and $R_1$ a positive constant. Suppose that $-1\leq R\leq R_1$ for all time $t$, and $u$ is a positive solution of the heat equation \eqref{ricciheateq}.
Then, for any $\delta\in[\frac{1}{2},1)$, we have
\begin{equation}\label{LYRF}
\delta\frac{|\nabla u|^2}{u^2}-\frac{\partial_t u}{u} \le \delta\frac{|\nabla u|^2}{u^2}-\frac{\partial_t u}{u}-\alpha R+\frac{\beta}{R+2}\leq
\frac{1}{t} \left(\frac{n}{2\delta}+\frac{4n\beta T}{\delta(1-\delta)} \right)
\end{equation}
for $t\in(0,T)$, where $\alpha=\frac{n}{2\delta(1-\delta)^2}$ and $\beta=\alpha (R_1+2)^2$.

\end{theorem}

Note that the curvature assumption is made only on the scalar curvature rather than on the Ricci or curvature
tensor. In this sense, this assumption is essentially optimal.  Under suitable assumptions, the
result in the theorem still holds when $\M$ is complete noncompact, or with a little weakening of the scalar curvature condition.
For the Ricci flow on a compact manifold $\M$, one can always rescale a solution so that the scalar curvature to be bounded from below by $-1$.

This theorem clearly implies a Harnack inequality for positive solutions of \eqref{ricciheateq}
if the scalar curvature is bounded.

 In the proof, the main tool is still the maximum principle applied on a differential inequality involving
Li-Yau type quantity. However, due to the Ricci flow, extra terms involving the Ricci curvature and
Hessian of the solution will appear. In order to proceed we need to create a new term with
the scalar curvature in the denominator. The key identity is the following.
Let
\be\label{defF}
F=-\Delta u+\delta\frac{|\nabla u|^2}{u}-\alpha Ru+\frac{\beta u}{R+C},
\ee and operator $\mathcal{L}=\Delta -\frac{\partial }{\partial t}$, where $\delta$, $\alpha$, and $\beta$ are arbitrary constants and $C$ is a constant so that $R+C>0$. Then
\be\lab{eqF}
\begin{aligned}
\mathcal{L}F&= \frac{1}{u}\left|u_{ij}-\frac{u_iu_j}{u}+uR_{ij}\right|^2+
\\
& \quad \frac{2\delta-1}{u}\left|u_{ij}-\frac{u_iu_j}{u}\right|^2+\frac{1}{(2\alpha-1)u}\left|(2\alpha-1)uR_{ij}+\frac{u_iu_j}{u}\right|^2\\
&\quad -\frac{1}{2\alpha-1}\frac{|\nabla u|^4}{u^3}+\frac{\alpha u}{R+C}\left|\nabla R-\frac{R+C}{u}\nabla u\right|^2\\
&+(\frac{\beta}{(R+C)^3}-\frac{\alpha}{R+C})u|\nabla R|^2
 -\frac{\alpha (R+C)|\nabla u|^2}{u}\\
 &+\frac{2\beta u|R_{ij}|^2}{(R+C)^2}+\frac{\beta u}{R+C}\left| \frac{\nabla R}{R+C}-\frac{\nabla u}{u}\right|^2
-\frac{\beta|\nabla u|^2}{u(R+C)}.
\end{aligned}
\ee Multiplying this by $t$ and using the maximum principle yields the claimed result.

\section{Li-Yau gradient bound under integral curvature condition}

In all of these log gradient bounds  in the static metric cases, the essential assumption is that the Ricci curvature or the corresponding Bakry-Emery Ricci curvature is bounded from below by a constant.
In many situations, it is highly desirable to weaken this assumption.
For example, in the case of K\"ahler Ricci flow in the Fano case, at each time slice,   one does not know if
the Ricci curvature is bounded. One only knows that $|Ric|$ is in $L^4$ by the work of
G. Tian and Z. L. Zhang \cite{TZz1} and $|Ric|^2$ is in a Kato class by G. Tian and Q.S. Zhang \cite{TZq2}.

By the work Petersen and Wei \cite{PW1} and \cite{PW2}, we know that the classical volume
comparison theorem can be generalized to the case with $L^p$ curvature bound with $p>n/2$.
If one also has Li-Yau bound, then one knows that Poincar\'e and Harnack inequality also hold.
These allows one to extend, to certain extent, the Cheeger-Colding theory to the case of integral Ricci bound.

 A few years ago,  Meng Zhu  and the author  proved Li-Yau bounds for positive solutions for both the fixed metric case  \eqref{heatequation} and  the Ricci flow case
\eqref{ricciheateq} under essentially optimal curvature conditions. The later was discussed in the last section.

For the fixed metric case, we will have two independent conditions
and two conclusions. The conditions are motivated by different problems such as studying
manifolds with integral Ricci curvature bound and the K\"ahler-Ricci flow. The conclusions
range from long time bound with necessarily worse constants, to short time
bound with better constants.

The theorem basically says:

(a).  $| Ric^-|$ in $L^p$, $p>n/2$  and volume noncollapsing
$\Rightarrow$
Li-Yau inequality $\Rightarrow$  classical parabolic Harnack inequality. i.e. if $u$ is a positive solution of the heat equation and $Q^- =B(x, r) \times [T_1, T_2]$ and $Q^+=B(x, r) \times [T_3, T_4]$ are two cylinders in the domain of $u$ such that $T_2<T_3$. Then
\[
\sup_{Q^-} u \le C \inf_{Q^+} u.
\]Here $C$ is a positive constant independent of $u$.

(b). Along each time slice of normalized K\"ahler Ricci flow, the Li-Yau inequality holds.

\begin{theorem} (\cite{ZZ1})\label{thmLY}
Let $(\M, g_{ij})$ be a compact $n$ dimensional Riemannian manifold, and $u$ a positive solution of \eqref{heatequation}.  Suppose either one of the following conditions holds.

(a)  $\int_\M |Ric^-|^p dy \equiv \sigma <\infty$ for some $p>\frac{n}{2}$, where $Ric^-$ denotes the nonpositive part of the Ricci curvature;  and the manifold is noncollapsed under scale 1, i.e., $|B(x, r)| \ge \rho r^n$ for $0<r\leq 1$ and some
 $\rho>0$;

(b)   $\sup_\M \int_\M |Ric^-|^2 d^{-(n-2)}(x,y)dy \equiv \sigma <\infty$ and the heat kernel of \eqref{heatequation}
satisfies the Gaussian upper bound (which holds automatically under (a)):
\be
\lab{GUP}
G(x, t; y, 0) \le \frac{\hat C(t)}{t^{n/2}} e^{ - \bar c d^2(x, y)/t},  \, t \in (0, \infty)
\ee for some positive constant $\bar c$ and positive increasing function $\hat C(t)$ which grows to infinity as $t\to \infty$.

Then, for any constant $\alpha\in(0,1)$, we have
\be\label{eqLY}
\a \underline{J}(t) \frac{ |\d u|^2}{u^2} - \frac{\p_t u}{u} \le \frac{n}{(2 - \delta)\a\underline{J}(t)} \frac{1}{t}
\ee
for $t\in (0, \infty)$, where $$\delta=\frac{2(1-\a)^2}{n+(1-\a)^2},$$ and
\be
\underline{J}(t)=\left\{
\al
&2^{-1/(5\delta^{-1}-1)}e^{-(5\delta^{-1}-1)^{\frac{n}{2p-n}}\left[4\sigma\hat C(t)^{1/p}\right]^{\frac{2p}{2p-n}}t},\ \textrm{under condition (a)};\\
&2^{-1/(10\delta^{-1}-2)}e^{-C_0(5\delta^{-1}-1)\sigma \hat C(t)t},\ \textrm{under condition (b)},
\eal\right.
\ee
with $C_0$ being a constant depending only on $n$, $p$ and $\rho$, and $\hat C(t)$ the increasing function on the right hand side of \eqref{GUP}.
\end{theorem}

In particular, for any $\beta\in(0,1)$, there is a $T_0=T_0(\beta,\sigma, p, n, \rho)$ such that
\be\label{eqLY2}
\beta\frac{ |\d u|^2}{u^2} - \frac{\p_t u}{u} \le \frac{n}{2\beta} \frac{1}{t}
\ee
for $t\in (0, T_0]$. Here $T_0 = c (1-\beta)^{4p/(2p-n)}$ and $c (1-\beta)^{4}$
under conditions (a) and (b), respectively; and $c$ is a positive constant depending only
on the parameters of conditions (a) and (b), i.e., $c=c(\sigma, p, n, \rho) $.

Are the conditions sharp?  It turns out that the $L^p$ condition on $|Ric^-|$ is almost necessary for solutions to be differentiable. For example one can not take $p=n/2$. The power $2$ on top of $|Ric^-|$ in condition (b) can be replaced by 1 and the proof is identical. It is also realized in \cite{Ro} and \cite{Car} that a variation of this integral condition on $|Ric^-|$, using the heat kernel and space time integration, instead of $1/d(x, y)^{n-2}$ and space integral, which is also called as the Kato condition, essentially implies heat kernel upper bound with a growing parameter.
 Moreover the volume noncollasping condition can be removed.

 For some constants $p,\ r>0$, define $\displaystyle k(p,r)=\sup_{x\in M}r^2\left(\oint_{B(x,r)}|Ric^-|^p dV\right)^{1/p}$, where $Ric^-$ denotes the negative part of the Ricci curvature tensor. We prove that for any $p>\frac{n}{2}$, when $k(p,1)$ is small enough, certain Li-Yau type gradient bound holds for the positive solutions of the heat equation on geodesic balls $B(O,r)$ in $\M$ with $0<r\leq 1$. Here the assumption that $k(p,1)$ being small allows the situation where the manifolds is collapsing.

\begin{theorem}(\cite{ZZ2})
Let $(\M^n, g_{ij})$ be a complete Riemannian manifold and $u$ a positive solution of the heat equation on $\M$.
For any $p>\frac{n}{2}$, there exists a constant $\kappa=\kappa(n,p)$ such that the following holds. If $k(p,1)\leq \kappa$, then for any point $O\in \M$ and constant $0<\alpha<1$, we have
\be\label{Li-Yau}
\a \underline{J}\frac{|\d u|^2}{u^2}-\frac{\partial_t u}{u}\leq \frac{n}{\a(2 - \delta)\underline{J}}\frac{1}{t}+\frac{C}{\a(2 - \delta)\underline{J}}\left[\frac{1}{\a(2-\delta)\underline{J}(1-\a)}+1\right],
\ee
in $B(O,\frac{1}{2})\times(0,\infty)$, where
\[
\underline{J}=\underline{J}(t)=2^{-\frac{1}{a-1}}\exp\left\{-2C\kappa\left(1+[2C(a-1)\kappa]^{\frac{n}{2p-n}}\right)t\right\},
\] $\delta=\frac{2(1-\a)^2}{n+(1-\a)^2}$, $a= \frac{5[n+(1-\alpha)^2]}{2(1-\alpha)^2}$ and $C=C(n,p)$ is a constant depending on $n$ and $p$.
\end{theorem}

The theorem basically says if $K(p, r)$ is small then Li-Yau inequality holds.
The smallness condition is necessary due to the example by Deane Yang \cite{Yan}, which is a dumb bell shaped manifold where $K(p, r)$ is not small.  There is no volume doubling property for such a manifold. But we know Li-Yau inequality implies, among several things, the volume doubling property.
See the discussion in a recent paper by X.Z. Dai, G.F. Wei and Z. L. Zhang \cite{DWZ}.

So the condition of this theorem is also essentially sharp.

In particular, on such manifolds  the $L^2$ Poincar\'e inequality holds. It has been an open question for a while if this is true for collapsing manifolds with integral Ricci bound.

One can summarize the result as:
$K(p, r)$  small  implies:

(a).  Harnack inequality for heat equation;

(b). volume doubling property and Poicar\'e inequality.

Notice that the only condition is on the Ricci curvature in terms of $K(p, r)$.  Previously by the work of C. Croke \cite{Cr},
one knows that the pointwise bound $Ricci \ge -k$ implies the isoperimetric inequality and hence the $L^2$ Sobolev inequality. In \cite{DWZ}, the authors gave a direct proof of the isoperimetric inequality under integral Ricci bound using relative volume comparison and an idea by Gromov \cite{Gro}. See also earlier result in \cite{Gal}.

From the work of Saloff-Coste \cite{Sa} and Grigoryan \cite{Gri}, we know that from the Harnack inequality, there are a number of other implications such as Gaussian upper bound of heat kernel and the $L^2$ Sobolev inequality.
 By the work of Petersen and Wei \cite{PW1}, we also know it implies a number of geometric and topological results such as finiteness of diffeomorphism type.

Let us outline the idea of the proof.
We will use the maximum principle on a quantity involving the solution and its gradient in the end.  However, unlike the case with Ricci
curvature bounded from below, the term involving the Ricci curvature can not be thrown away.
The idea is to construct an axillary function by solving a nonlinear equation to cancel the term
involving the Ricci curvature.  The function is now referred to as the $J$ function.

Let $J=J(x, t)$ be a smooth positive  function and $\a \in (0, 1)$ be a parameter.
Write
\be
\lab{defforQ}
Q \equiv \a J \frac{ |\d u|^2}{u^2} - \frac{\p_t u}{u}.
\ee and $f=\ln u$. After some computation one finds that
\be
\lab{ineqtQ}
\al
&(\Delta  - \p_t) ( t Q )
+ 2 \d \ln u \d (t Q)
\ge  \a t (2 J - \delta J) \frac{1}{n}  \left( | \d  f |^2 - \p_t f \right)^2 \\
&\quad  +
\a \left[\Delta J - 2 V  J
-
5 \delta^{-1} \frac{| \d J|^2}{J}
 -\p_t J \right]  t |\d f|^2
 - \delta \a  t J |\d f|^4 -Q,
\eal
\ee where we have written $| Ric^-|=V$ and $\delta>0$ is another small parameter.

For any given parameter $\delta>0$ such that $5\delta^{-1}>1$,
the problem
\be
\lab{eqforJ}
\begin{cases}
\Delta J - 2 V  J
-
5 \delta^{-1}\frac{| \d J|^2}{J}
 -\p_t J =0, \quad \text{on} \quad {\M} \times (0, \infty);\\
J(\cdot, 0) = 1,
\end{cases}
\ee
has a unique solution for $t\in[0,\infty)$, which satisfies
\be
\underline{J}(t)\leq J(x,t)\leq 1,
\ee
where
\be
\underline{J}(t)=\left\{
\al
&2^{-1/(5\delta^{-1}-1)}e^{-(5\delta^{-1}-1)^{\frac{n}{2p-n}}\left[4\sigma\hat C(t)^{1/p}\right]^{\frac{2p}{2p-n}}t},\ \textrm{under condition (a)};\\
&2^{-1/(10\delta^{-1}-2)}e^{-C_0(5\delta^{-1}-1)\sigma \hat C(t)t},\ \textrm{under condition (b)},
\eal\right.
\ee
with $C_0$ being a constant depending only on $n$, $p$ and $\rho$, and $\hat C(t)$ the increasing function in \eqref{GUP}. The term involving the Ricci curvature on \eqref{ineqtQ} is cancelled. The maximum principle then implies the result in the theorem.

\section{Li-Yau gradient bound with sharp constants}

Except for the case where Ricci is nonnegative and Perelman's gradient bounds for the fundamental solution, these Li-Yau type bounds are not sharp.
The difference between
\eqref{Li-Yau.op} and \eqref{Li-Yau>1} is more than a formality since the former is actually a sharp Log Laplace estimate because it says
\[
- \Delta \ln u = \frac{|\nabla u|^2}{u^2} - \frac{u_t}{u} \le \frac{n}{2t}.
\]

The Li-Yau bound \eqref{Li-Yau>1} when Ricci curvature changes sign is not sharp.

   The bound \eqref{Li-Yau>1} was later improved for small time by Hamilton in \cite{Ha93}, where he proved under the same assumptions as above that
\begin{equation}\label{Hamilton}
\frac{|\nabla u|^2}{u^2} - e^{2Kt}\frac{u_t}{u} \leq e^{4Kt}\frac{n}{2t}.
\end{equation}
See also the work by Davies \cite{Dav} and Yau himself [Yau94].

For the original Li-Yau inequality, in over 30 years, several sharpening of the bounds have been obtained with $\alpha$ replaced by several functions $\alpha=\alpha(t)>1$ but not equal to $1$. We have mentioned some of these work in Section 2. A well known open question  ( in Chow-Lu-Ni's book, \cite{LX}, \cite{BBG}  and several other places) asks if a sharp bound can be reached. In a recent short paper, we show that for all complete compact manifolds one can take $\alpha=1$. Thus a sharp bound, up to computable constants, is found in the compact case. This result also seems to sharpen Theorem 1.4 in \cite{LY:1} for compact manifolds with convex boundaries.

In the noncompact case one can not take $\alpha=1$ even for the hyperbolic space. An example is also given, which shows that there does not exist an optimal function $\alpha=\alpha(t, K)$ for all noncompact manifolds with negative Ricci lower bound $-K$, giving a negative answer to the open question in the noncompact case.

\begin{theorem}(\cite{Z21})
\label{main thm}
Let $(\M, g_{ij})$ be a complete, compact $n$ dimensional Riemannian manifold and $u$ a positive solution of the heat equation on $\M \times (0, \infty)$, i.e.,
\be\label{HE}
(\Delta - \partial_t)u=0.
\ee Let $diam_\M$ be the diameter of $\M$ and suppose the Ricci curvature is bounded from below by a non-positive constant  $-K$, i.e. $R_{ij} \ge -K g_{ij}$. Then there exist  dimensional constants $C_1>0$ and $C_2>0$ such that
\be
\lab{sly}
\al
&-t \Delta \ln u = t \left(\frac{|\nabla u|^2}{u^2} - \frac{\partial_t u}{u} \right)\\
&\le \left( \frac{n}{2} + \sqrt{2n K (1+K t)(1+t)} \,  diam_\M + \sqrt{ K (1+K t) (C_1  + C_2 K ) t} \right).
\eal
\ee
\end{theorem}

The constants $C_1$ and $C_2$ arise from the standard volume comparison or/and heat kernel upper and lower bound which can be efficiently estimated. For compact manifolds, the main utility of the result is for small or finite time. For example, when we choose $C_1$ and $C_2$ to be independent of the volume comparison theorem, the bound actually gives rise to a Laplace and hence volume comparison theorem as $t \to 0$.  For large time, a solution converges to a constant. As $t \to \infty$, it is expected that the bound will be of order $o(t)$ instead of $O(t)$ as of now. But we will not pursue the improvement this time. It would still be interesting to find a sharp form of the RHS of \eqref{sly}.

For  the noncompact case, we give a negative answer on the existence of a sharp function which depends on time $t$ and Ricci lower bound only and for which the Li-Yau bound holds. First, in the noncompact case one can not prove a bound like \eqref{sly} with a finite right hand side for a fixed $t$ even for the hyperbolic space. For example in the 3 dimensional hyperbolic space with the standard metric, the heat kernel (c.f. \cite{Dav} e.g.) is:
 \[
 G(x, t, y) = \frac{1}{(4 \pi t)^{3/2}} \frac{  d(x, y)}{\sinh d(x, y)} e^{-t- \frac{d^2(x, y)}{4t}}.
 \]For $r=d(x, y)$,
 \[
 \al
 &-t \Delta \ln G(x, t, y) =\\
 &- t \sinh^{-2} r \frac{\partial}{\partial r} \left(\sinh^2 r \frac{\partial}{\partial r} \right) \ln \left(\frac{1}{(4 \pi t)^{3/2}} \frac{r}{\sinh r} e^{-t- \frac{r^2}{4t}}\right)
 \eal
 \]is of the order $d(x, y)$ when $d(x, y) \to \infty$.

Therefore one can not take $\alpha=1$ for all noncompact manifolds with Ricci bounded from below.
 Second, let $\M_0$ be any compact manifold of dimension $n$ and $\M=\M_0 \times \R^1$ be the product manifold of $\M_0$ with $\R^1$. Then the heat kernel on $\M$ is given by $G=G_0(x, t, y) G_1(z, t, w)$ where $G_0$ and $G_1$ are the heat kernel on $\M_0$ and $\R^1$ respectively with $x, y \in \M_0$ and $z, w \in \R^1$ and $t>0$. Fixing $y$ and $w$, we compute
 \[
 \al
 &-t \Delta \ln G = -t (\Delta_0 + \partial^2_z) \ln G = -t\Delta_0 \ln G_0 -  t \partial^2_z \ln G_1 \le \\
 & \left( \frac{n+1}{2} + \sqrt{2n K (1+K t)(1+t)} \,  diam_{\M_0} + \sqrt{ K (1+K t) (C_1  + C_2 K ) t} \right).
 \eal
 \]

Here we just used the above Theorem  on the compact manifold $\M_0$ with $-K$ being the Ricci lower bound and $\Delta_0$ being the Laplace-Beltrami operator on $\M_0$. So by a result of C.J. Yu and F.F. Zhao \cite{YZ} or \cite{LWH+} ,  for the product manifold $\M$, we can take $\alpha=1$ for all positive solutions, which is optimal.

 These two examples show that there does not exist a single optimal function of time $t$ and $K$ only, which works for all noncompact manifolds with Ricci curvature bounded from below by a negative constant $-K$.

 The reason is that for each negative constant $-K$, the optimal function for the above manifold $\M_0 \times \R^1$ is $\alpha=1$. But for the hyperbolic space with $Ricci = - K g$, the optimal function is worse than $\alpha=1$.

The idea of the proof is the following. First we prove the sharp bound for heat kernels $G$ instead of general positive solutions. This involves an iterated integral estimate on $-\Delta \ln G$ which satisfies a nonlinear heat type equation and Hamilton's estimate of $|\nabla \ln G|$ by $t \ln A/G$, taking advantage of the finiteness of the diameter of $\M$ and lower bound of $G$. Here $A$ is a constant which dominates the heat kernel in a suitable time interval. Note that we do not follow the usual way of using the maximum principle. Next we use crucially the result in Yu-Zhao \cite{YZ} which states that a Li-Yau type bound on the heat kernel implies the same bound on all positive solutions.

Recently an improvement is made by Xingyu Song, Ling Wu, and Meng Zhu on the parameters,  resulting in some refinement and extension of the above result,
 Here is a sample in the Ricci flow case. See also \cite{SW} for another extension.

\begin{theorem} (\cite{SWZ})
\label{th1.5}
	Let $M^n$ be an n-dimensional closed Riemannian manifold and $(M^n,g(x,t))_{t\in [0,2)}$ a solution to the Ricci flow  $\frac{\partial}{\partial t}g(x,t)=-2 Ric(x,t)$, $ (x,t)\in M\times [0,2)$. Assume that for some constant $K\ge0$, $-\frac{K}{2-t}g(x,t)\le Ric(x,t)\le \frac{K}{2-t}g(x,t)$ and $\nu[g_0,4]\ge -nK $, $  (x,t)\in M\times [0,2)$. Let $u$ be a smooth positive solution of  the heat equation $\left(\Delta_t-\frac{\partial}{\partial t}\right)u(x,t)=0$ on $M\times [0,2)$ and $d_t(M)$ the diameter of $(M^n,g(x,t))$. Then there exists a constant $c_9$ depending  on $n,K$ such that
	\begin{equation}
		t\left(\frac{|\nabla u|^2}{u^2}-\frac{\partial_t u}{u}\right) \le \frac{n}{2}+ \frac{nKt}{2-t}+\sqrt{\frac{Kc_9}{2-t}\left(t+d_0^2(M)t+d_0^2(M)\right)}\nonumber
	\end{equation}
	for all $(x,t)\in M\times [0,2)$.
\end{theorem}

This sharpens the result of  \cite{BCP} where the coefficient in front of $|\nabla u|^2/u^2$ is not 1 and Ricci curvature is assumed to be bounded.  Instead, even some singularity of the Ricci curvature is allowed.

Some possible problems to work on are.
\begin{itemize}

\item  Characterize noncompact manifolds satisfying the sharp Li-Yau bound?

\item sharp Li-Yau bound for the forward heat equation under Ricci flow?

\item  sharp Li-Yau bound in large time for the forward heat equation in noncompact case?

\item  Hamilton's entropy bound Theorem \ref{thhament} under integral Ricci condition?

\end{itemize}

\section{Matrix log gradient bounds}

As mentioned, R. Hamilton extended the Li-Yau estimate to the log matrix form.
On a general compact $n$-manifold $\M$, he proved the following

\begin{theorem} (\cite{Ha93})
\lab{thham93}
 For any positive solution $u$ to \eqref{heatequation}, there exist constants $B$ and $C$ depending only on the geometry of $\M$ (in particular the diameter, the volume, and the curvature and covariant derivative of the Ricci curvature) such that if $t^{n/2}u\leq B$ then
\begin{equation}\label{eq Hamilton matrix general case}
\n_i \n_j \ln u  +\frac{1}{2t}g_{ij} +C \left(1+\ln\left(\frac{B}{t^{n/2}u}\right)\right)g_{ij} \geq 0.
\end{equation}

Moreover $C=0$ if $\M$ has nonnegative sectional curvature and parallel Ricci curvature and the inequality holds for all positive solution $u$.

\end{theorem}

Joint with Xiaolong Li, we found an improved version of the above theorem by Hamilton.

\begin{theorem}( \cite{LZ})
\label{thm improve Hamilton}
Let $(\M,g)$ be a closed Riemannian $n$-manifold and let $u:\M \times [0,T] \to \R$ be a positive solution to the heat equation \eqref{heatequation}. Suppose that the sectional curvatures of $M$ are bounded by $K$ and $|\n Ric|\leq L$, for some $K,L>0$. Then
\begin{eqnarray}
 \n_i \n_j \log u  +
 \left(\frac{1}{2t}+(2n-1)K +\frac{\sqrt{3}}{2}L^{\frac{2}{3}} +\frac{1}{2t}\gamma(t,n,K,L)\right)g_{ij} \geq 0
\end{eqnarray}
for all $(x,t)\in \M \times (0,T)$, where $\gamma(t,n,K,L)$ is
\begin{eqnarray*}
  && \sqrt{nKt(2+(n-1)Kt)}  +\sqrt{C_3(K+L^{\frac{2}{3}})t(1+Kt)(1+K+Kt)}\\
&& + \left(2K(2+(n-1)Kt)+\frac{3}{2}L^{\frac{2}{3}}(1+(n-1)Kt) \right) Diam,
\end{eqnarray*}
$C_3>0$ depends only on the dimension $n$, and $Diam$ denotes the diameter of $(\M,g)$.
\end{theorem}

Notice that Hamilton's original inequality \eqref{eq Hamilton matrix general case} has the term $C\log \left(\frac{B}{t^{\frac n 2}u}\right)$, where $B$ and $C$ depend on the geometry of the manifold and $B$ is greater than $t^{\frac n 2}u$, which itself is an additional assumption. The constant $C$ is equal to zero only when $\M$ has nonnegative sectional curvature and parallel Ricci curvature. Otherwise, for this $\log$ term, we do not have any definite control on the order $q$ of $-t^{-q}$ coming out of this term, for general positive solutions, making this lower bound less practical. In Theorem \ref{thm improve Hamilton}, we manage to replace this term with a $C/t$ term, with $C$ depending only on $K$, $L$, and $Diam$, which is of the correct order for $t$.

We also establish  a similar result for a general compact Ricci flow.

\begin{theorem} (\cite{LZ})
\label{thm matrix Harnack heat equation general case}
Let $(\M,g(t))$, $t\in [0,T]$, be a compact solution to the Ricci flow.
Let $u:\M \times [0,T] \to \R$ be a positive solution to the heat equation
\eqref{ricciheateq}.
Suppose the sectional curvatures of $(\M,g(t))$ are bounded by $K$ for some $K>0$.
Then
\begin{equation}\label{eq matrix Harnack general case}
        \n_i \n_j \ln u + \left(\frac{1}{2 t} + \frac{1}{t} \beta(t, n, K) \right) g_{ij} \ge 0,
\end{equation}
where
$$\beta(t, n, K)= 4 \sqrt{n K t}  +C_2 (K+1)t +C_1 \sqrt{K} diam.$$ Here $C_1>0$ is a numerical constant, $C_2>0$ depends only on the dimension and the non-collapsing constant $v_0 = \inf \{|B(x, 1, g(0))|_{g(0)}: x \in M \}$, and
$$diam:= \sup_{t\in [0,T]} diam(M,g(t)).$$
\end{theorem}

Similar upper bounds for $ Hess \, \ln u$ can be found in \cite{HZ}. The idea in that paper of using local extension of eigen-vectors into vector fields is also instrumental in the proof of the preceding two theorems. This method allows us to avoid Hamilton's global tensor maximum principle.

\begin{theorem}(\cite{HZ})
\label{ThmHessfixed} Let $\M$ be a Riemannian $n$-manifold
with a metric $g$.

$\operatorname{(a)}$ Suppose $u$ is a solution of
$$\partial_{t} u- \Delta u =0\quad \text{in} \ \M \times (0, T].$$
Assume $0 < u \leq A$. Then, $$ t ( \n_i \n_j u) \leq
u(5+Bt)\left( 1+ \log \frac{A}{u}\right)\quad \text{in} \ \M \times (0, T],$$
where $B$
is a nonnegative constant depending only on the
$L^\infty$-norms of curvature tensors, $n$
and the gradient of Ricci curvatures.

$\operatorname{(b)}$
Suppose $u$ is a
solution of
$$\partial_{t} u- \Delta u = 0\quad\text{in }Q_{R,T}(x_0,t_0)=B(x_0, R) \times [t_0-T, t_0].$$
Assume $0 < u \leq A$. Then,
$$( \n_i \n_j u)\leq C u \left( \frac{1}{T}
+ \frac{1}{R^{2}} + B\right) \left(1+\log \frac{A}{u}\right)^{2}
\quad \text{in} \ Q_{\frac{R}{2}, \frac{T}{2}}(x_0,t_0).$$ where
$C$ is a universal constant and $B$ is a nonnegative constant
depending only on the $L^\infty$-norms of curvature tensors, $n$ and
the gradient of Ricci curvatures in $Q_{R,T}(x_0,t_0)$.
\end{theorem}

Notice that $ \n_i \n_j \ln u = \n_i \n_j  u /u  - \n_i   u \n_j u /u^2.$ So the above theorem together with Theorem \ref{thhament} and Theorem \ref{THSZ1} yield the respective upper bound for $ \n_i \n_j \ln u$.

Hamilton's Theorem \ref{thham93} is sharp when the sectional curvature is nonnegative and the Ricci curvature is parallel. The equality is reached by the heat kernel in the Euclidean space. It is noticed by Poon \cite{Po} that this sharp matrix estimate is connected with the unique continuation property for the heat equation on such manifolds.

In the paper \cite{LZ}, we  extended Hamilton's sharp matrix estimate
\eqref{thham93}, $C=0$ case for static metrics to the Ricci flow case. Our estimate does not require the parallel Ricci curvature condition and thus could be more applicable.

\begin{theorem}\label{thm matrix Harnack heat equation}
Let $(\M^n,g(t))$, $t\in [0,T]$, be a complete solution to the Ricci flow.
Let $u:\M^n \times [0,T] \to \R$ be a positive solution to the heat equation
\begin{equation}\label{heat equation}
u_t -\Delta_{g(t)}u=0.
\end{equation}
Suppose that $(\M^n,g(t))$ has nonnegative sectional curvature and $Ric \leq \kappa g$ for some constant $\kappa>0$. Then
\begin{equation}\label{matrix estimate K geq 0}
        \n_i \n_j
        \log u + \frac{\kappa }{1-e^{-2\kappa t}} g_{ij} \geq 0,
\end{equation}
    for all $(x,t)\in M\times (0,T)$.
\end{theorem}

In the same paper, this theorem is used to prove unique continuation property (UCP) of the conjugate heat equation on those manifolds using monotonicity of a modified parabolic frequency function with the heat kernel of the forward heat equation as weight. Let us recall that a solution $u$ satisfies (UCP) if vanishing at infinity order at one space time point implies vanishing everywhere. This result is in part inspired by the results in the papers \cite{Lin90}, \cite{Po}, \cite{CM} and \cite{BK}. See also \cite{LW19} and \cite{LLX22}. It is well known that (UCP) requires very sharp parameters in the monotonicity formula, which is provided by Theorem \ref{thm matrix Harnack heat equation}. An interesting question is whether the theorem still holds without the upper bound on the Ricci curvature. Also the parameters in the matrix Harnack estimates are likely not sharp except we have nonnegative sectional curvature and parallel Ricci curvature in the static case. It is desirable to find sharp ones.

\medskip
\hspace{-.5cm}{\bf Acknowledgements}
The author gratefully acknowledges the support of Simons'
Foundation grant 710364 and Prof. Xiaolong Li for helpful suggestions. He also wishes to thank Professors Jiayu Li, Gang Tian, Zhenlei Zhang and  Xiaohua Zhu for the invitation to write the paper and to the annual geometric analysis conference.


\begin{thebibliography}{99}
{\small
%\bibitem[AbGr]{AbGr} Abresch, Uwe; Gromoll, Detlef, {\it On complete manifolds with nonnegative Ricci curvature}, J. Amer. Math. Soc. 3 (1990), no. 2, 355-374.
%\bibitem[And]{Anderson} Anderson, Michael T., {\it A survey of Einstein metrics on 4-manifolds}, Handbook of geometric analysis, No. 3, 1-39, Adv. Lect. Math. (ALM), 14, Int. Press, Somerville, MA, 2010, http://arxiv.org/abs/0810.4830.


\bibitem[An]{An}  Anderson, M., {\it On the topology of complete manifolds of nonnegative Ricci curvature.} Topology 29 (1990),
no. 1, 41-55.

\bibitem[ArBe]{ArBe} Aronson, Donald G.; B\'enilan, Philippe, {\it R\'egularit\'e des solutions de l'\'equation des milieux poreux dans RN}, C. R. Acad. Sci. Paris S\'er. A-B 288 (1979), no. 2, A103-A105.

\bibitem[Ber]{Ber} S. N. Bernstein, {\it Sur la g\'en\'eralisation du probl\'eme de Dirichlet. 1st, 2nd part, } Math. Ann. 62 (1906), 253-271; Math. Ann. 69 (1910), 82-136.

\bibitem[BCP]{BCP} Bailesteanu, Mihai; Cao, Xiaodong; Pulemotov, Artem Gradient estimates for the heat equation under the Ricci flow. J. Funct. Anal. 258 (2010), no. 10, 3517-3542.

\bibitem[BBG]{BBG} Bakry, Dominique; Bolley, Francois; Gentil Ivan, {\it The Li-Yau inequality and applications under a curvature-dimension condition},  arXiv:1412.5165, 2014.

\bibitem[BL]{BL} Bakry, Dominique and Ledoux, Michel,
{\it A logarithmic Sobolev form of the Li-Yau parabolic inequality},
  Rev. Mat. Iberoamericana, Volume 22, Number 2 (2006), 683-702.


\bibitem[BNS]{BNS} Elia Brue, Aaron Naber, Daniele Semola,
{\it Fundamental Groups and the Milnor Conjecture},
arXiv:2303.15347

\bibitem[BZ]{BZ} Bamler, Richard; Zhang, Qi S., {\it Heat kernel and curvature bounds in Ricci flows with bounded scalar curvature}, Advances in Mathematics 319, 396-450, 2017
%\bibitem [BCG]{BCG:1} Barlow, Martin; Coulhon, Thierry; Grigor'yan, Alexander, {\it Manifolds and graphs with slow heat kernel decay}, Invent. Math. 144 (2001), no. 3, 609-649.

\bibitem[Bre]{Bre} Brendle, Simon, {\it A generalization of Hamilton's differential Harnack inequality for the Ricci flow}, J. Differential Geom. 82 (2009), no. 1, 207-227.

\bibitem[BK]{BK} Julius Baldauf and Dain Kim. {\it Parabolic frequency on Ricci flows,} Int. Math. Res.
Not. IMRN, to appear, arXiv:2201.05505v2, (2022).

\bibitem[Cao]{Cao}Cao, Huai-Dong, {\it On Harnack's inequalities for the K\"ahler-Ricci flow}, Invent. Math. 109 (1992), no. 2, 247-263.

\bibitem[CaNi]{CaNi} Cao, Huai-Dong; Ni, Lei, {\it Matrix Li-Yau-Hamilton estimates for the heat equation on K\"ahler manifolds}, Math. Ann. 331 (2005), no. 4, 795-807.


\bibitem[CTZ]{CTZ} Cao, Huai-Dong; Tian, Gang; Zhu, Xiaohua, {\it K\"ahler-Ricci solitons on compact complex manifolds with $C_1(M)>0$}. Geom. Funct. Anal. 15 (2005), no. 3, 697-719.


\bibitem[Cx]{Cx} Cao, Xiaodong, {\it Differential Harnack estimates for backward heat equations with potentials under the Ricci flow.}  J. Funct. Anal. 255 (2008), no. 4, 1024-1038.

%\bibitem[CaHa]{CaHa} Cao, Xiaodong; Hamilton, Richard S., {\it Differential Harnack estimates for time-dependent heat equations with potentials}, Geom. Funct. Anal. 19 (2009), no. 4, 989-1000.

\bibitem [CH]{CH:1} Cao, Xiaodong; Hamilton, Richard S., {\it Differential Harnack estimates for time-dependent heat equations with potentials},  Geom. Funct. Anal. 19 (2009), no. 4, 989-1000.

%\bibitem[CZ]{CZ} Cao, Xiaodong; Zhang, Qi S. {\it The Conjugate Heat Equation and Ancient Solutions of the Ricci Flow}, Adv. Math., Vol. 228 (2011), no. 5, 2891-2919.
%\bibitem[CG]{CG} Cheeger, Jeff; Gromoll, Detlef, {\it The splitting theorem for manifolds of nonnegative Ricci curvature}, J. Differential Geometry 6 (1971/72), 119-128.

%\bibitem[ChCo]{ChCo} Cheeger, Jeff; Colding, Tobias H., {\it Lower bounds on Ricci curvature and the almost rigidity of warped products}, Ann. of Math. (2) 144 (1996), no. 1, 189-237.

\bibitem[Car]{Car} Carron, Gilles, {\it  Geometric inequalities for manifolds with Ricci curvature in the Kato class}, Ann. Inst. Fourier (Grenoble) 69 (2019), no. 7, 3095-3167.

\bibitem[Chow]{Chow} Chow, Bennett, {\it The Ricci flow on the 2-sphere}, J. Differential Geom. 33 (1991), no. 2, 325-334.

\bibitem [ChH]{ChH:1} Chow, Bennett; Hamilton, Richard S., {\it Constrained and linear Harnack inequalities for parabolic equations}, Invent. Math. 129 (1997), no. 2, 213-238.

%\bibitem[CN]{Cheeger-Naber-codim-4} Cheeger, Jeff; Naber, Aaron, {\it Regularity of Einstein Manifolds and the Codimension 4 Conjecture}, http://arxiv.org/abs/1406.6534 (2014).

%\bibitem[CT]{CH:2} Cao, Xiaodong; Tran, Hung, {\it Mean Value Inequalities and Conditions to Extend Ricci Flow}, http://arxiv.org/abs/1303.4492v1 (2013).



%\bibitem [CW1]{CW:1} Chen, Xiuxiong; Wang, Bing, {\it Space of Ricci flows I}, Comm. Pure Appl. Math. 65 (2012), no. 10, 1399-1457, http://arxiv.org/abs/0902.1545

\bibitem [CW13]{CW13} Chen, Xiuxiong; Wang, Bing, {\it On the conditions to extend Ricci flow(III)},
Int. Math. Res. Not. IMRN 2013, no. 10, 2349-2367.

%\bibitem [CW3]{CW:3} Chen, Xiuxiong; Wang, Bing, {\it Space of Ricci flows (II)}, http://arxiv.org/abs/1405.6797



\bibitem[Chowetc]{Chowetc} Chow, Bennett; Lu, Peng; Ni, Lei, {\it Hamilton's Ricci flow}. Graduate Studies in Mathematics, 77. American Mathematical Society, Providence, RI; Science Press, New York, 2006. xxxvi+608 pp.

\bibitem[CoNa]{CoNa} Colding, Tobias Holck; Naber, Aaron, {\it Sharp H\"older continuity of tangent cones for spaces with a lower Ricci curvature bound and applications}, Ann. of Math. (2) 176 (2012), no. 2, 1173-1229.

\bibitem[Cr]{Cr} C. Croke, {\it Some isoperimetric inequalities and eigenvalue estimates}, Ann. Sci. Ecole Norm. Sup. , 13 (1980) pp. 419-435.

\bibitem[CM]{CM} Tobias Holck Colding and William P. Minicozzi, II. {\it  Parabolic frequency on manifolds.}
Int. Math. Res. Not. IMRN, (15): 11878-11890, (2022).

\bibitem[CY]{CY} S. Y. Cheng and S. T. Yau,  {\it Differential equations on Riemannian manifolds
and their geometric applications}, Comm. Pure Appl. Math. 28:3 (1975), 333-354.

%\bibitem[GH]{GH:2} Grigor'yan, Alexander; Hu, Jiaxin, {\it Off-diagonal upper estimates for the heat kernel of the Dirichlet forms on metric spaces}, Invent. Math. 174 (2008), no. 1, 81-126.

\bibitem[Dav]{Dav} E.  B.  Davies, {\it
Heat  kernels  and  spectral  theory}, volume  92  of
Cambridge  Tracts  in Mathematics,   Cambridge University Press, Cambridge, 1989

\bibitem[DWZ]{DWZ} Xianzhe Dai, Guofang Wei, Zhenlei Zhang, {\it Local Sobolev Constant Estimate for Integral Ricci Curvature Bounds}, Advances in Mathematics 325 (2018) 1-33.

 \bibitem[Gal]{Gal} S. Gallot, {\it Isoperimetric Inequalities Based on Integral Norms of Ricci Curvature}, Asterisque,
vol. 157-158, Soc. Math. de France, 1988, pp. 191-216.

\bibitem[GM]{GM} N. Garofalo and A. Mondino, {\it  Li-Yau and Harnack type inequalities in
$RCD^*(K;N)$ metric measure spaces.} Nonlinear Anal., 95: 721-734,  2014.

\bibitem[Gri]{Gri} A. A. Grigor'yan, {\it The heat equation on noncompact Riemannian manifolds}, Mat. Sb., 182:1 (1991), 55-87; Math. USSR-Sb., 72:1 (1992), 47-77

%\bibitem[Gri]{Gr1997} Grigor'yan, Alexander, {\it Gaussian upper bounds for the heat kernel on arbitrary manifolds}, J. Differential Geom. 45 (1997), no. 1, 33-52.

%\bibitem[Gr]{Gr:2} Grigor'yan, Alexander, {\it Heat kernel and analysis on manifolds}, AMS/IP Studies in Advanced Mathematics, 47. American Mathematical Society, Providence, RI; International Press, Boston, MA, 2009. xviii+482 pp.

%\bibitem [Ha1]{Ha:0}Hamilton, Richard S., {\it Three-manifolds with positive Ricci curvature}, J. Differential Geom. 17 (1982), no. 2, 255-306.

%\bibitem [Ha2]{Ha:1} Hamilton, Richard S., {\it The formation of singularities in the Ricci flow}, Surveys in differential geometry, Vol. II (Cambridge, MA, 1993), 7-136, Int. Press, Cambridge, MA, 1995.

\bibitem[Gro]{Gro} M. Gromov, {\it Paul Levy's Isoperimetric Inequality}, Appendix C in Metric Structures for Riemannian
and Non-Riemannian Spaces, Progr. Math., vol. 152, Birkhauser, 2001.

\bibitem[Ha93]{Ha93} Hamilton, Richard S., {\it A matrix Harnack estimate for the heat equation}, Comm. Anal. Geom. 1 (1993), no. 1, 113-126.

\bibitem[Ha4]{Ham1988} Hamilton, Richard S., {\it The Ricci flow on surfaces}, Mathematics and general relativity (Santa Cruz, CA, 1986), 237-262, Contemp. Math., 71, Amer. Math. Soc., Providence, RI, 1988.

\bibitem[Ha5]{Ham1993jdg}  Hamilton, Richard S., {\it The Harnack estimate for the Ricci flow}, J. Differential Geom. 37 (1993), no. 1, 225-243.

%\bibitem[HN]{HN:1} Hein, Hans-Joachim; Naber, Aaron, {\it New logarithmic Sobolev inequalities and an $\varepsilon$-regularity theorem for the Ricci flow}, Comm. Pure Appl. Math. 67 (2014), no. 9, 1543-1561.

\bibitem[HY]{HY} Bobo Hua, Wenhao Yang, {\it Liouville theorems for ancient solutions of subexponential growth to the heat equation on graphs}, arXiv:2309.17250.

\bibitem[HZ]{HZ} Qing Han, Qi S Zhang, {\it An upper bound for Hessian matrices of positive solutions of heat equations}, Journal of Geometric Analysis 26 (2016), 715-749.

\bibitem[KuZh]{KuZh} Kuang, Shilong; Zhang, Qi S., {\it A gradient estimate for all positive solutions of the conjugate heat equation under Ricci flow}, J. Funct. Anal. 255 (2008), no. 4, 1008-1023.

\bibitem[Li]{Li} Li, Peter, {\it Large time behavior of the heat equation on complete manifolds with nonnegative Ricci curvature.}, Ann. of
Math. (2) 124 (1986), no. 1, 1-21

\bibitem [Lib] {Lib} Li, Peter,
{\it Geometric analysis.} Cambridge Studies in Advanced Mathematics, 134. Cambridge University Press, Cambridge, 2012. x+406 pp.

\bibitem[Liu]{Liu} Liu, Gang, {\it 3-manifolds with nonnegative Ricci curvature}. Invent. Math. 193 (2013), no. 2, 367-375

\bibitem[LLX22]{LLX22}    Chuanhuan Li, Yi Li, and Kairui Xu, {\it Parabolic frequency monotonicity on Ricci
flow and Ricci-harmonic flow with bounded curvatures.} J. Geom. Anal., to appear,
arXiv:2205.07702, (2022).

\bibitem[LW19]{LW19} Xiaolong Li and Kui Wang, {\it Parabolic frequency monotonicity on compact manifolds},
Calc. Var. Partial Differential Equations, 58(6):Paper No. 189, 18, (2019).

 \bibitem[LX]{LX} J. Li and X. Xu.  {\it Differential Harnack inequalities on Riemannian manifolds I: linear
heat equation.} Adv. Math., 226(5):4456-4491, 2011.

\bibitem[LZ]{LZ} Xiaolong Li, Qi S. Zhang, {\it Matrix Li-Yau-Hamilton estimates under Ricci Flow and parabolic frequency}, Calc. Var. PDE, 2024;  arXiv:2306.10143.




\bibitem[Lij]{Lij} J.Y. Li, {\it Gradient estimate for the heat kernel of a complete Riemannian manifold and its applications}, J. Functional Analysis 97, 293-310 (1991)

%\bibitem[L]{L:1} Li, Peter, {\it Geometric analysis}, Cambridge Studies in Advanced Mathematics, 134. Cambridge University Press, Cambridge, 2012.

\bibitem[Lin90]{Lin90} Fang-Hua Lin. {\it A uniqueness theorem for parabolic equations.} Comm. Pure Appl.
Math., 43(1):127-136, (1990).

\bibitem[LinZ]{LinZ} Fanghua Lin and Q. S. Zhang. {\it On ancient solutions of the heat equation.} Comm. Pure Appl.
Math., 72(9):2006-2028, 2019.

\bibitem[LWH+]{LWH+} Zhen Lei, Zhiqiang Wang, Bobo Hua, Peng Qu, and Genggeng Huang, {\it Lectures on
mathematical and physical equations}. Fudan Lecture Notes, in Chinese, unpublished

\bibitem[LT]{LT} Li, Peter; Tian, Gang, {\it On the heat kernel of the Bergmann metric on algebraic varieties}, J. Amer. Math. Soc., 8 (1995), 857-877.

\bibitem [LY]{LY:1} Li, Peter; Yau, Shing-Tung, {\it On the parabolic kernel of the Schr\"odinger operator.}  Acta Math. 156 (1986), no. 3-4, 153-201.

\bibitem[Ma]{Ma} T. Mabuchi, {\it Heat kernel estimates and the Green functions on Multiplier
hermitian manifolds}, Tohoku Math. J. 54 (2002), 259-275.

\bibitem[Mo]{Mo} Modica, Luciano, {\it  A gradient bound and a Liouville theorem for nonlinear Poisson equations}, Comm. Pure Appl. Math. 38 (1985), no. 5, 679-684.

\bibitem[Mos]{Mos} Sunra Mosconi. {\it Liouville theorems for ancient caloric functions via optimal growth conditions}.
Proc. Amer. Math. Soc., 149(2):897-906, 2021.

\bibitem[Ni]{Ni2006} Ni, Lei, {\it A note on Perelman's LYH-type inequality}, Comm. Anal. Geom. 14 (2006), no. 5, 883-905.

\bibitem[Pa]{Pa} Pan, Jiayin, {\it A proof of Milnor conjecture in dimension 3}. J. Reine Angew. Math. 758 (2020), 253-260.

\bibitem[P1]{P:1} Perelman, Grisha, {\it The entropy formula for the Ricci flow and its
geometric applications}, http://arxiv.org/abs/math/0211159 (2002).

%\bibitem[P2]{P:2} Perelman, Grisha, {\it Ricci flow with surgery on three-manifolds}, http://arxiv.org/abs/math/0303109 (2003).

\bibitem[Po]{Po} Chi-Cheung Poon, {\it Unique continuation for parabolic equations. Comm. Partial Differential Equations}, 21(3-4):521-539, 1996.

\bibitem[PW1]{PW1} Petersen, Peter; Wei, Guofang, {\it Relative volume comparison with integral curvature bounds}, Geom. Funct. Anal. 7 (1997), no. 6, 1031-1045.

\bibitem[PW2]{PW2} Petersen, Peter; Wei, Guofang, {\it Analysis and geometry on manifolds with integral Ricci curvature bounds. II}, Trans. Amer. Math. Soc. 353 (2001), no. 2, 457-478.


\bibitem[QZZ]{QZZ}
  Z. Qian, H.-C. Zhang, and X.-P. Zhu.  {\it Sharp spectral gap and Li-Yau's estimate on
Alexandrov spaces}.
Math. Z., 273(3-4):1175-1195, 2013

%\bibitem[Se]{Sesum-Ricci}Sesum, Natasa, {\it Curvature tensor under the Ricci flow}, Amer. J. Math. 127 (2005), no. 6, 1315-1324.

%\bibitem[Sh]{Sh:1} Shi, Wan-Xiong, {\it Deforming the metric on complete Riemannian manifolds}, J. Differential Geom. 30 (1989), no. 1, 223-301.

%\bibitem[Si]{Si:1} Simon, Miles, {\it  Ricci flow of almost non-negatively curved three manifolds}, J. Reine Angew. Math. 630 (2009), 177-217.

%\bibitem[ST]{ST:1} Sesum, Natasa; Tian, Gang, {\it Bounding scalar curvature and diameter along the K\"ahler Ricci flow (after Perelman)}, J. Inst. Math. Jussieu 7 (2008), no. 3, 575-587.

%\bibitem[TW]{Tian-Wang} Tian, Gang; Wang, Bing, {\it On the structure of almost Einstein manifolds}, http://arxiv.org/abs/1202.2912 (2012).

\bibitem[Ro]{Ro} Rose, Christian, {\it Li-Yau gradient estimate for compact manifolds with negative part of Ricci curvature in the Kato class}, Ann. Global Anal. Geom. 55 (2019), no. 3, 443-449.

\bibitem[Sa]{Sa} Saloff-Coste, L. {\it A note on Poincar\'e, Sobolev, and Harnack inequalities},
International Mathematics Research Notices, Volume 1992, Issue 2, 1992, Pages 27-38.

\bibitem[So]{So} Sormani, C., {\it Nonnegative Ricci curvature, small linear diameter growth and finite generation of fundamental
groups}. J. Differential Geom. 54 (2000), no. 3, 547-559.

\bibitem[SW]{SW}  Xingyu Song, Ling Wu, {\it Li-Yau gradient estimates on closed manifolds under bakry-emery ricci curvature conditions}, arXiv:2204.12851.


\bibitem[SWZ]{SWZ} Xingyu Song, Ling Wu, Meng Zhu, {\it A direct approach to sharp Li-Yau Estimates on closed manifolds with negative Ricci lower bound}, arXiv:2307.03879.

\bibitem[SZ]{SZ:1} Souplet, Philippe; Zhang, Qi S. {\it
Sharp gradient estimate and Yau's Liouville theorem for the heat equation on noncompact manifolds.}
 Bull. London Math. Soc. 38 (2006), no. 6, 1045-1053.

 \bibitem[TZq1]{TZq1} Tian, Gang; Zhang, Qi S.  {\it Isoperimetric inequality under K\"ahler Ricci flow}.
 Amer. J. Math. 136 (2014), no. 5, 1155-1173.

\bibitem[TZq2]{TZq2} Tian, Gang; Zhang, Qi S., {\it A compactness result for Fano manifolds and K\"ahler Ricci flows},  Math. Ann.  362  (2015),  no. 3-4, 965-999.

\bibitem[TZz]{TZz1} Tian, Gang; Zhang, Zhenlei,
{\it Regularity of  K\"ahler Ricci flows on Fano manifolds}, Acta Math. 216(1): 127-176 (2016).

%\bibitem[Wa]{Wa:1} Wang, Bing, {\it On the conditions to extend Ricci flow (II)}, Int. Math. Res. Not. IMRN 2012, no. 14, 3192-3223.

%\bibitem[Ye]{Ye:1} Rugang Ye, {\it The logarithmic Sobolev inequality along the Ricci flow}, http://arxiv.org/abs/0707.2424 (2007).
\bibitem[Wan]{Wan}  F.-Y. Wang. {\it Gradient and Harnack inequalities on noncompact manifolds with boundary},
Pacific J. Math., 245(1):185-200, 2010.

\bibitem[WanJ]{WanJ} Wang, Jiaping, {\it Global heat kernel estimates.} Pacific J. Math. 178 (1997), no. 2, 377-398.

\bibitem[Ya]{Ya} S. T. Yau, {\it Harmonic functions on complete Riemannian manifolds}, Comm. Pure Appl.
Math. 28 (1975), 201-228

\bibitem[Yan]{Yan} D. Yang, {\it Convergence of Riemannian manifolds with integral bounds on curvature}. I, Ann. Sci. Ec.
Norm. Super. 25 (1992) 77-105.

\bibitem[YZ]{YZ} Chengjie Yu, Feifei Zhao, {\it Li-Yau multiplier set and optimal Li-Yau gradient estimate on hyperbolic spaces}, arXiv:1807.05709.

\bibitem[Z06]{Z06:1} Zhang, Qi S., {\it Some gradient estimates for the heat equation on domains and for an equation by Perelman}, Int. Math. Res. Not., 39 pages Art. ID 92314, 39, 2006.

\bibitem[Z12]{Z12} Zhang, Qi S., {\it Bounds on volume growth of geodesic balls under Ricci flow}, Math. Res. Letters, Volume 19, Issue 1 (2012), pp. 245-253

\bibitem[Z21]{Z21} Qi S. Zhang. {\it A sharp Li-Yau gradient bound on compact manifolds.} arXiv:2110.08933 v2,
2021

\bibitem[ZZ1]{ZZ1} Q.S. Zhang and M. Zhu, {\it Li-Yau gradient bounds on compact manifolds under nearly optimal curvature conditions,}
Journal of Functional Analysis 275 (2), 478-515, 2018.

\bibitem[ZZ2]{ZZ2} Q.S. Zhang and M. Zhu, {\it Li-Yau gradient bound for collapsing manifolds under integral curvature condition},
Proceedings of the American Mathematical Society 145 (7), 3117-3126, 2017.}

%\bibitem[Z07]{Z07:1} Zhang, Qi S., {\it A uniform Sobolev inequality under Ricci flow}, IMRN 2007, ibidi Erratum, Addendum.

%\bibitem[Z11]{Z11:1} Zhang, Qi S., {\it Bounds on volume growth of geodesic balls under Ricci flow}, Math. Res. Lett. 19 (2012), no. 1, 245-253, http://arxiv.org/abs/1107.4262

%\bibitem[ZZh]{ZhouZhang-scal} Zhang, Zhou, {\it Scalar curvature behavior for finite-time singularity of K\"ahler-Ricci flow}, Michigan Math. J. 59 (2010), no. 2, 419-433.
\end{thebibliography}
\end{document}